\documentclass[11pt]{amsart}

\usepackage{tasks}
\usepackage{hyperref}
\usepackage{blindtext}
\usepackage{setspace}

\usepackage{graphicx}
 
\usepackage{subcaption}
\usepackage[utf8]{inputenc}
\usepackage[export]{adjustbox}
\usepackage{wrapfig}

\usepackage{marginnote}
\usepackage{datetime}
\usepackage{color}
\usepackage{amsrefs}
\usepackage{amsmath}
\usepackage{mathtools}
\usepackage{graphicx}  
\usepackage[top=3.5cm, bottom=2.5cm, outer=2.5cm, inner=2.5cm]{geometry}
\usepackage{hyperref}
\hypersetup{linkcolor=red,
citecolor=black
}
\usepackage{amsmath}

\newtheorem{thm}{Theorem}[section]
\newtheorem{lem}[thm]{Lemma}

\theoremstyle{definition}

\newtheorem{Alg}[thm]{Algorithm}

\usepackage{hyperref}
\usepackage{hyperref}
\hypersetup{
	colorlinks = true,
	linkcolor = {red},
	citecolor = {blue}
}

\theoremstyle{remark}
\newtheorem{re}[thm]{Remark}

\newcommand{\RNum}[1]{\uppercase\expandafter{\romannumeral #1\relax}}

\numberwithin{equation}{section}

\usepackage[utf8]{inputenc}
\usepackage[english]{babel}
\usepackage{paralist}
\usepackage{amsthm}
\usepackage{marginnote}

\theoremstyle{definition}

\theoremstyle{remark}

 \usepackage{bm,amsmath}

\newcommand{\tran}{^{\mathstrut\scriptscriptstyle\top}} 
 
\newcommand{\RN}[1]{%
  \textup{\uppercase\expandafter{\romannumeral#1}}%
}

\newcommand{\Rey}{\mathcal{R}e }

\newcommand{\grad}{\nabla}
\newcommand{\bfzero}{\mathbf{0}}
\newcommand{\bff}{\mathbf{f}}
\newcommand{\bfu}{\mathbf{u}}
\newcommand{\bfv}{\mathbf{v}}
\newcommand{\bfw}{\mathbf{w}}
\newcommand{\bfX}{\mathbf{X}}

\newcommand{\bfV}{\mathbf{V}}
\newcommand{\bfH}{\mathbf{H}}
\newcommand{\Cinv}{\tilde C}
\def\bfF{\mathbf{F}}

\allowdisplaybreaks[4]

\usepackage[utf8]{inputenc}

\title{Data assimilation with higher order finite element interpolants}

\author{Michael S. Jolly}
\author{Ali Pakzad}
\thanks{The research of M.~Jolly was supported in part by NSF grant DMS-1818754. The authors acknowledge the Lilly Endowment, Inc., through its support for the Indiana University Pervasive Technology Institute, which provided supercomputing resources used for this research \url{https://doi.org/10.5967/K8G44NGB}.
 }
 
 \date{\today}

\begin{document}
\begin{abstract}
The efficacy of a nudging data assimilation algorithm using higher order finite element interpolating operators is studied.  Numerical experiments are presented for the 2D Navier-Stokes equations in two cases: shear flow in an annulus and a forced flow in a disk with an off-center cavity.  In both cases second order interpolation of coarse-grain data is shown to outperform first order interpolation. Convergence of the nudged solution to that of a direct numerical reference solution is proved.  The analysis points to a trade-off in the estimates  for higher order interpolating operators.
\end{abstract}
\maketitle

\section{Introduction}
In rough terms, data assimilation refers to a class of methodologies which observational data is combined with a model in order to  improve the accuracy in forecasts. These  techniques have been used for a long time in weather modeling, climate science, and hydrological
and environmental forecasting  \cite{K03}. There are a variety of data assimilation techniques, whereby actual measured quantities over time are incorporated in system models.  Best known perhaps is the Kalman filter, which for discrete in time, linear systems give exact probabilistic predictions consistent with uncertanties in the both the data and the model.  This Bayesian approach is adapted to nonlinear systems in extended Kalman filters, though they are no longer exact.   Variational methods, 3DVar and 4DVar also track such uncertainties for discrete systems.  More on these methods can be found in the books \cites{Asch-Data2016, Evensen-Data2009, Law-AMathematical2015}. 

Continuous, deterministic data assimilation injects an interpolant of the data, assumed to be known for all time in some interval, directly into a differential equation.  
Charney, Halem and Jastrow proposed inserting low Fourier mode data into the advective term of the Navier-Stokes equations (NSE), \cites{Charney-Use1969, Henshaw-Numerical2003}. 
That approach can be interpreted through the result of Foias and Prodi \cite{Foias-Sur1967} showing that there are a finite number of determining modes for the 2D NSE.  Quantities such as volume elements and nodal values have also been shown to be determining \cites{Cockburn-Estimating1997, Constantin-Determining1985,Jones-Upper1993}. While these quantities are more readily measured in practice, inserting them into the advective term is problematic.   The technique of nudging 
avoids the need for derivatives by using this data in a linear feedback term.
For a model given by
\begin{align}\label{origODE} 
\frac{d \bfu}{dt} =\bfF(\bfu)
\end{align}
whose solution is known over a coarse grid of resolution $H$ as $I_H \mathbf{u}(t)$, nudging is done through the auxiliary system
\begin{align}\label{nudgedODE} 
\frac{d \bfv}{dt} =\bfF(\bfv) -\mu I_H(\bfv-\bfu)\;, \quad \bfv(0)=0\;.
\end{align} 
There is literature devoted to a form of this technique for synchronization of chaotic dynamical systems.  See \cites{Auroux-ANudging2008, Pecora-Synchronization2015} for a more complete history dating back to \cite{Hoke-TheInitialization1976}, and \cites{Pawar-Long2020,Vidard-Determination2003} for comparisons with Kalman filtering.  

For \eqref{origODE} given by the 2D NSE, Azouani, Olson and Titi  showed that for large enough $\mu$ and small enough $H$, $\|\bfv-\bfu\| \to 0$ at an exponential rate \cite{AOT14}.   Such analysis has since been completed for a variety of PDEs modeling physical phenomena  \cites{Biswas-Continuous2018,Jolly-Adata2017,Jolly-Continuous2019,Jolly-Determining2017,Markowich-Continuous2019,Pei-Continuous2019, MR3930850, ZRSI19}, in some instances, using data in only a subset of system variables \cites{Farhat-Continuous2015, Farhat-Continuous2017, Farhat-Data2020, Farhat-Abridged2016, Farhat-Data2016} and more recently, over a  subdomain of the physical spatial domain \cite{Biswas-Data2021}.   
Though deterministic, the nudging approach has been studied with respect to error in the observed data \cites{Bessaih-Continuous2015, Jolly-Continuous2019}.  Computational experiments have demonstrated that nudging is effective using data that is much more coarse than the rigorous estimates require  \cites{Altaf-Downscaling2017,Farhat-Assimilation2018,Gesho-Acomputational2016,Hudson-Numerical2019}, supported in some cases by analysis of certain numerical schemes
\cites{Foias-ADiscrete2016, Mondaini-Uniform2018}. 


Nudging has also been studied for use with finite element method (FEM). Uniform in time error estimates in the semi-discrete case are made in \cite{Garcia-Archilla-Uniform2020}.   Fully discrete FEM nudging schemes are shown to be well posed and stable in \cite{LRZ19}.  Both works analyze the error between the reference solution of the 2D NSE and the finite dimensional numerical scheme, and test the efficacy numerically in cases where the exact reference solution is specified over the square $[0,1]^2$, and the corresponding body force added.   A numerical test of 2D channel flow past a cylinder is also done in \cite{LRZ19}, in which case the reference solution is taken to be the result of a direct FEM numerical simulation over a fine mesh.
Both works also focus on the case where interpolating operator satisfies
\begin{align}\label{zeroorder}
\|I_H \bfu -\bfu\|_{L^2} \le c \,H\, \|\nabla \bfu\|_{L^2}\;,
\end{align}
which is achieved by taking constant values over each triangle in a FEM discretization.

In this paper we present numerical evidence to show that higher order interpolation can achieve faster synchronization by nudging.  This is suggested by the bound 
\begin{equation}\label{Ineq:InterpolateIntro}
 \| \phi  - I^k_H\phi\|_m\leq C_{k+1,m}\,  H^{k+1-m} \| \, \phi\|_{ k+1},
 \end{equation}
for $k$-degree polynomial over each triangle \cite{BR08}. Two flows are tested: a shear flow in an annulus, and one with a body force in a disk with an off-center obstacle, both satisfying Dirichlet boundary conditions.   For the shear flow we demonstrate that using a quadratic interpolating polynomial ($k=2$), nudging can succeed
when the data is too coarse for using a linear polynomial ($k=1$).  When the data is finer, both synchronize to within machine precision, but the quadratic polynomial does so in less than one-third the time.  Similarly, in the body force case, for data with a certain resolution $H$, the speed-up is about a factor of 2 when a quadratic interpolating polynomial is used.
As we do not have exact reference solutions for these flows, the errors are in terms of the difference between the nudged solution and that of a direct numerical simulation (DNS) on a fine mesh, triangles of diameter $h\ll H$. Synchronization with the true solution of the NSE by nudging with higher order interpolation is analyzed in \cite{Biswas-Higher2021}.  Here, we include a proof in the semi-discrete case of convergence of the nudged solution to the DNS solution for interpolating polynomials of any order.   While this limited analysis is not sensitive enough to show an advantage in using higher degree interpolation, it does indicate a trade-off between the higher power in \eqref{Ineq:InterpolateIntro} and the price one pays to close the estimates with an inverse inequality \eqref{Ineq:Inverse}.  The computational results provide evidence that the higher power wins.

\section{Notation}
Let $\Omega$ be an open, bounded region in $\mathbb{R}^2$ with a Lipschitz continuous boundary.    Define the velocity space $\bfX$ as 

$$ \bfX := \bfH_0^1 (\Omega) =  \{\bfu \in \bfH^1(\Omega): \bfu = 0 \hspace{0.2cm}  \mbox{on} \hspace{0.2cm}   \partial \Omega \},$$
and the pressure space $Q$  as follows

$${Q}(\Omega):= L_0^2(\Omega) = \{q \in L^2(\Omega): \hspace{3pt}\int_{\Omega} q\  dx= 0\}.$$
The closed subspace of divergence free functions is given by
$$ \bfV :=   \{\bfu \in \bfX: \, (\nabla \cdot \bfu , q ) = 0, \, \, \forall \ q \in Q \}.$$
We denote the (explicitly skew symmetrized)  trilinear form $b : \bfX \times \bfX  \times \bfX  \rightarrow \mathbb{R}$ as 

$$b (\bfu, \bfv, \bfw) := \frac{1}{2} (\bfu \cdot \nabla \bfv, \bfw) - \frac{1}{2} (\bfu \cdot \nabla \bfw, \bfv),$$
which has the following property
$$b (\bfu, \bfv, \bfw)  = - b (\bfu, \bfw, \bfv),$$
and therefore
$$b (\bfu, \bfv, \bfv) = 0, \hspace{1cm} \forall \bfu, \bfv \in \bfX.$$

  To discretize the equations, consider a regular mesh $\Omega_h \in \Omega$  with maximum triangle diameter length $h$. Let two finite-dimensional spaces $\bfX^h \subset \bfX$ and $Q^h \subset Q $ be finite element velocity and pressure spaces corresponding to an admissible triangulation of $\Omega$. We further assume that $ \bfX^h $ and $Q^h$  satisfy the following  discrete inf-sup condition (the condition for div-stability)  
\begin{equation}\label{LBBH}
\inf_{q^h \in \mathbb{Q}^h} \sup_{\bfv^h \in \bfX^h} \frac{(q^h, \nabla \cdot \bfv^h)}{\|\nabla \bfv^h\|_0 \, \|q^h\|_0} \geq \beta ^h >0,
\end{equation}
where $\beta ^h >0$ uniformly in $h$ as $h \rightarrow 0$.

In most finite element discretizations of the NSE and related systems, the divergence-free constraint $\nabla\cdot \bfu^h = 0$ is only weakly enforced. What holds instead of the pointwise constraint is that a numerical solution $\bfu^h$ in $\bfX^h$ satisfies $\bfu^h \in  \bfV^h$, where $\bfV^h \subset \bfX^h$ is constructed as
$$\bfV^h:= \{ \bfv^h \in \bfX^h: (q^h, \grad \cdot \bfv^h)=0, \,  \,\,\forall  \, q^h \in Q^h \}.$$ 

\subsection{Triangular Finite Elements} Consider the most common finite element mesh; triangular.   We first give here two examples of conforming finite element basis functions that are used in our computation,  for more  details see \cite{BR08}.

\begin{enumerate}
\item   \textit{$P_1: C^0-$piecewise linear on triangles}  $(\text{dim}\,P_1 =3)$. On each element the basis function is of the form 
$$\phi(x,y) = c_0 + c_1 x + c_2 y,$$
and thus the values of $\phi$ are uniquely determined once $\phi$ is specified at the three vertices. 

\item   \textit{$P_2: C^0-$piecewise quadratic on triangles} $(\text{dim}\, P_2 =6)$. On each triangle, the basis functions are full quadratic polynomials
$$\phi(x,y) = c_0 + c_1 x + c_2 y + c_3 x^2 + c_4 y^2 + c_5 xy,$$
and thus six nodes are needed per element to determine $\phi$. The standard nodes are chosen to be the vertices and midpoints of the edges of the triangle.




\end{enumerate}

 One way to construct a velocity-pressure element spaces which satisfy (\ref{LBBH}) is by considering 

$$\bfX^h = \{\bfv^h \in \bfH_0^1 (\Omega) : \bfv^h |_{\bigtriangleup} \in P_k(\bigtriangleup) \, , \forall \,  \bigtriangleup \text{in the mesh}\},$$
$$Q^h = \{q^h \in L_0^2(\Omega)\cap C^0(\Omega) : q^h |_{\bigtriangleup} \in P_{k-1}(\bigtriangleup) \, , \forall \,  \bigtriangleup \text{in the mesh}\}.$$
The choice $k=2$, known as the  Taylor-Hood elements \cite{TaylorHood}, is one commonly used choice of velocity-pressure finite element spaces which satisfy the  discrete inf-sup condition.  In fact for any $k\geq 2$ the above choice of velocity-pressure element spaces  $\bfX^h  , Q^h$ satisfies (\ref{LBBH}) \cites{G89, J16, L08}.

Throughout this manuscript,  the $L^2(\Omega)$ and   $H^1(\Omega)$ norms  will be denoted by $\| \cdot \|_0$ and, $\| \cdot \|_1$  respectively, and the $L^2(\Omega)$  inner product  is given by $(\cdot \, , \,  \cdot)$.   For $k\geq 2$, the cell-wise definition of the  norms in Sobolev spaces $H^k(\Omega)$ will be considered since in these cases finite element functions do not possess the regularity for the global norm to be well defined.  Therefore for $k \geq 2$ without loss of generality we denote 
$$ \| \cdot \|_{k} = \sum_{\bigtriangleup \in \text{Mesh}}\,\| \cdot \|_{H^k(\bigtriangleup)}.$$
Let $\Omega_H$ denote a coarse finite element mesh which is refined by successively joining midpoints with line segments, ultimately producing the finest mesh $\Omega_h$, so $h \ll  H$. 
 The spacing $h$ corresponds to our Direct Numerical Simulation (DNS),  whose solution, $\bfu^h$, plays the role of the reference  solution, with which we seek to synchronize. The spacing $H$  corresponds in practice to points where the true solution is observed and the data is collected.  We process the observables by a $k$-degree  interpolating polynomial  $I^k_H : \bfH^{k+1} (\Omega)\rightarrow L^2(\Omega)$, satisfying the following approximation property:
 \begin{equation}\label{Ineq:Intrpolate}
 \| \phi  - I^k_H\phi\|_m\leq C_{k+1,m}\,  H^{k+1-m} \| \, \phi\|_{ k+1},
 \end{equation}
 for every $\phi \in \bfH^{k+1}$,  where $C_{k+1,m}$ is independent of $\phi$ and the mesh, see \cite{BR08}.  The case of linear interpolation, $k=1$, was considered for finite element discretizations in \cites{Garcia-Archilla-Uniform2020,LRZ19} and is consistent with the original assumption in \cite{AOT14}.   In this paper we  study the relative efficacy of using higher order interpolants in nudging. 
 Each mesh is made sufficiently regular to satisfy the inverse inequality 
  \begin{equation}\label{Ineq:Inverse}
\|\phi^h\|_{ k} \leq \Cinv_{k,m} \, h^{m-k}\, \|\phi^h\|_m,
 \end{equation}
 for every $\phi^h \in  \bfX^h,$ where $\Cinv_{k,m}$ is a constant independent of  $\phi^h$ and the mesh \cite{BR08}.  The inequalities \eqref{Ineq:Intrpolate}, \eqref{Ineq:Inverse} point to the trade-off in estimating the efficacy of nudging with higher order interpolants.  While the higher power in $H$ is beneficial, the inverse relation with $h$ works against us.  Some analysis balancing these two effects is given in section \ref{SecAnalysis}.

\section{Equations}

We consider the  incompressible Navier-Stokes equations (NSE) with Dirichlet boundary condition
\begin{equation}
 \label{NSE}
 \begin{split}
\bfu_t+   (\bfu \cdot \nabla) \bfu   - \nu \Delta \bfu +  \nabla p  & = \bff,\\ 
\nabla \cdot \bfu  &= 0,
\end{split}
\end{equation}
in the physical domain $\Omega \subset \mathbb{R}^2$. Here $\bfu$ is the velocity, $\bff = \bff(x, t)$ is the known body force, $p$ is the pressure, and $\nu$ is the kinematic viscosity.   Let $I_H$ be be an interpolation operator satisfying \eqref{Ineq:Intrpolate}, then corresponding data assimilation algorithm is given by the system
\begin{equation}
 \label{DANSE}
 \begin{split}
\bfv_t+   (\bfv \cdot \nabla) \bfv   - \nu \Delta \bfv +  \nabla q  & = \bff - \mu\, I_H(\bfv -\bfu) ,\\ 
\nabla \cdot \bfv  &= 0,
\end{split}
\end{equation}
where the initial value of $\bfv_0 $ is arbitrary.

 We consider two flows, one with a body force $\bff$, and one without. 
In the case of body force, let $\lambda_1$ be the smallest eigenvalue of the Stokes operator,  our complexity parameter is the Grashof number
 \begin{equation}\label{Grashof}
 G = \frac{1}{\nu^2 \, \lambda_1} \|\bff \|_0,
 \end{equation}
 otherwise it is the Reynolds number given  as $\Rey \sim 1/\nu$.

 We use an IMEX (implicit-explicit) scheme as the temporal discretization to avoid the solution of a nonlinear problem at each time step. In short, in this scheme,  at time $t = t_{n+1}$ the nonlinear term $\bfv_{n+1} \cdot \nabla \bfv_{n+1} $ is replaced by $\bfv_{n} \cdot \nabla \bfv_{n+1} $,  where $ \bfv_{n}  $ is obtained from already computed solution at time $t = t_{n}$. The spatial discretization is the finite element method and we discretize in space via the  Taylor-Hood element pair, although  it can be extended also to any LBB-stable pair easily.

\begin{Alg}\label{Algorithm}
  Given body force $\bff \in L^{\infty} ((0, \infty), L^2(\Omega))$,  initial condition $\bfv_0 \in \bfV_h$ , reference solution $\bfu^h_{n+1}$, and  $(\bfv^h_n, q_n^h) \in \bfX^h \times Q^h$,  compute  $(\bfv_{n+1}^h, q_{n+1}^h) \in \bfX^h \times Q^h$ satisfying
\begin{multline}\label{IMEX-FEM}
(\frac{\bfv_{n+1}^{h} - \bfv_{n}^{h}}{\Delta t},\Theta ^{h}) + b(\bfv^h_n,\bfv_{n+1}^{h},\Theta^{h}) + \nu \, (\nabla \bfv_{n+1}^{h},\nabla \Theta^{h})  - (q_{n+1}^{h}, \nabla \cdot \Theta^{h})
\\ =  (\bff_{n+1},\Theta^{h})  - \mu \, (I_H(\bfv^h_{n+1} - \bfu_{n+1}), \Theta^h),
\end{multline}
\begin{equation}\label{incomp}
(\nabla \cdot \bfv_{n+1}^{h},r^{h}) = 0,
\end{equation}
for all $(\Theta^{h},  r^{h}) \in  \bfX^h \times Q^h.$
\end{Alg}
The approximating solution $\bfv^h$ with an arbitrary initial condition is computed using Algorithm \ref{Algorithm}, while $I_H (\bfu_{n+1})$  represents our observations of the system at a coarse spatial resolution at the time step  $ t_{n+1}$. Note that at the time $t_{n+1}$ we have already obtained the relevant  data about the true solution, therefore $I_H(\bfu_{n+1})$ is interpreted   to be the most recent data. The well-posedness  of the above algorithm (\ref{Algorithm}) with constant interpolation, $k=0$,  is studied in \cite{LRZ19}.  

\section{Computational study \RNum{1}; Shear Flow} 
One of the classic problems in experimental fluid dynamics is shear flow, where the boundary condition is tangential.  As in the experiment of G.I. Taylor and M.M.A. Couette, we consider a 2D slice of flow between rotating cylinders \cite{F95}.   The domain $\Omega$ is an annulus 
$$\Omega = \{ (x,y) \in \mathbb{R}^2: x^2 + y^2 \leq 1^2 \hspace{0.3cm} \text{and}  \hspace{0.3cm} x^2 + y^2 \geq (0.1)^2 \}.$$
The flow is driven by the rotational force at the outer circle in an absence of body force $\bff = 0$, with no-slip boundary conditions  imposed on the inner circle. The Reynolds number is taken to be  $ \Rey = 600$.    In all experiments, the algorithms are implemented by using public domain finite element software FreeFEM++ \cite{FreeFEM}, and was run  on the Carbonate supercomputer at Indiana University.

\subsection{Spatial Discretization} 
\begin{figure}[h!]
  \centering
  \begin{subfigure}[b]{0.3\linewidth}
    \includegraphics[width=\linewidth]{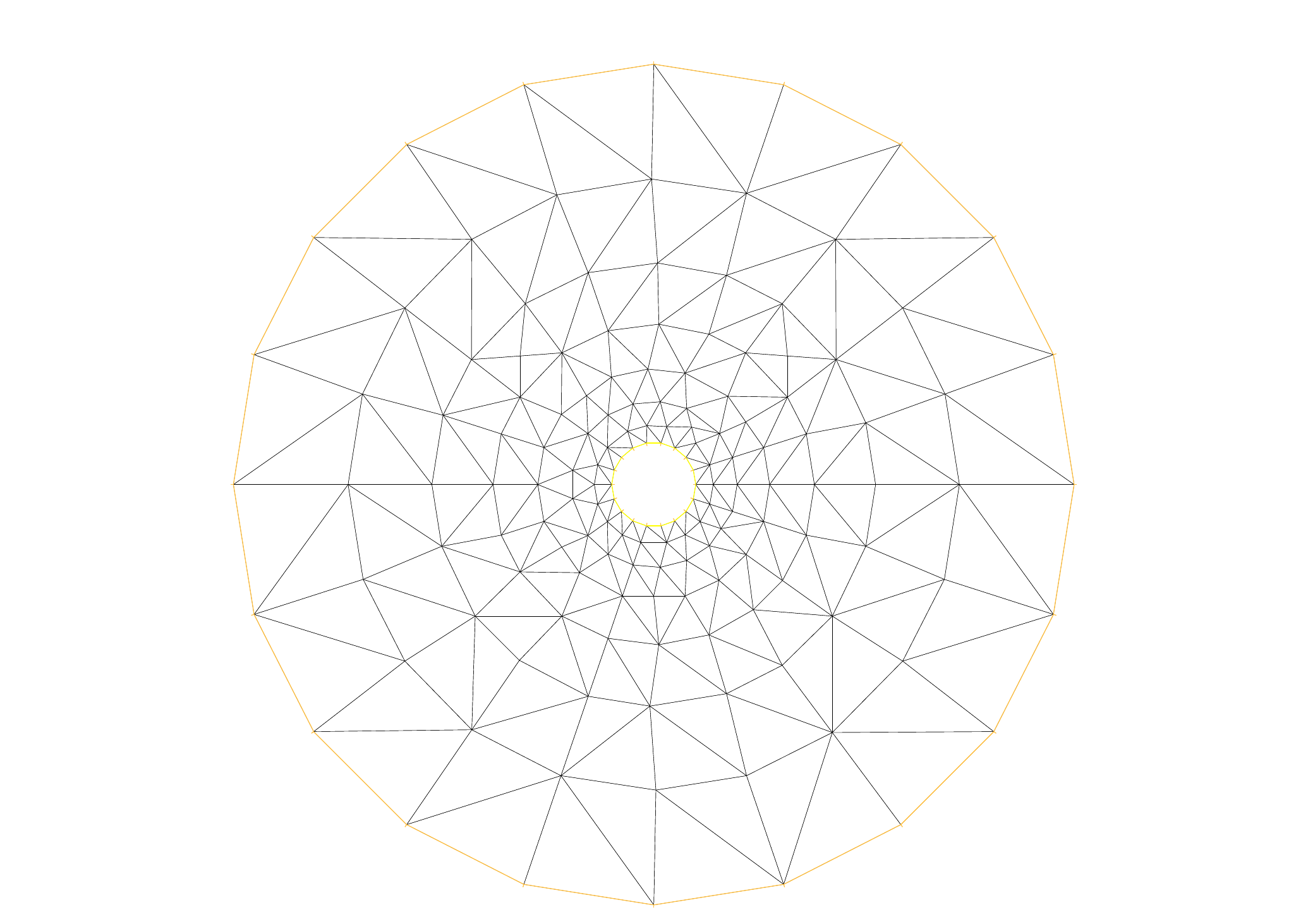}
    \caption{Mesh Level 1; Coarsest}
  \end{subfigure}
  \begin{subfigure}[b]{0.3\linewidth}
    \includegraphics[width=\linewidth]{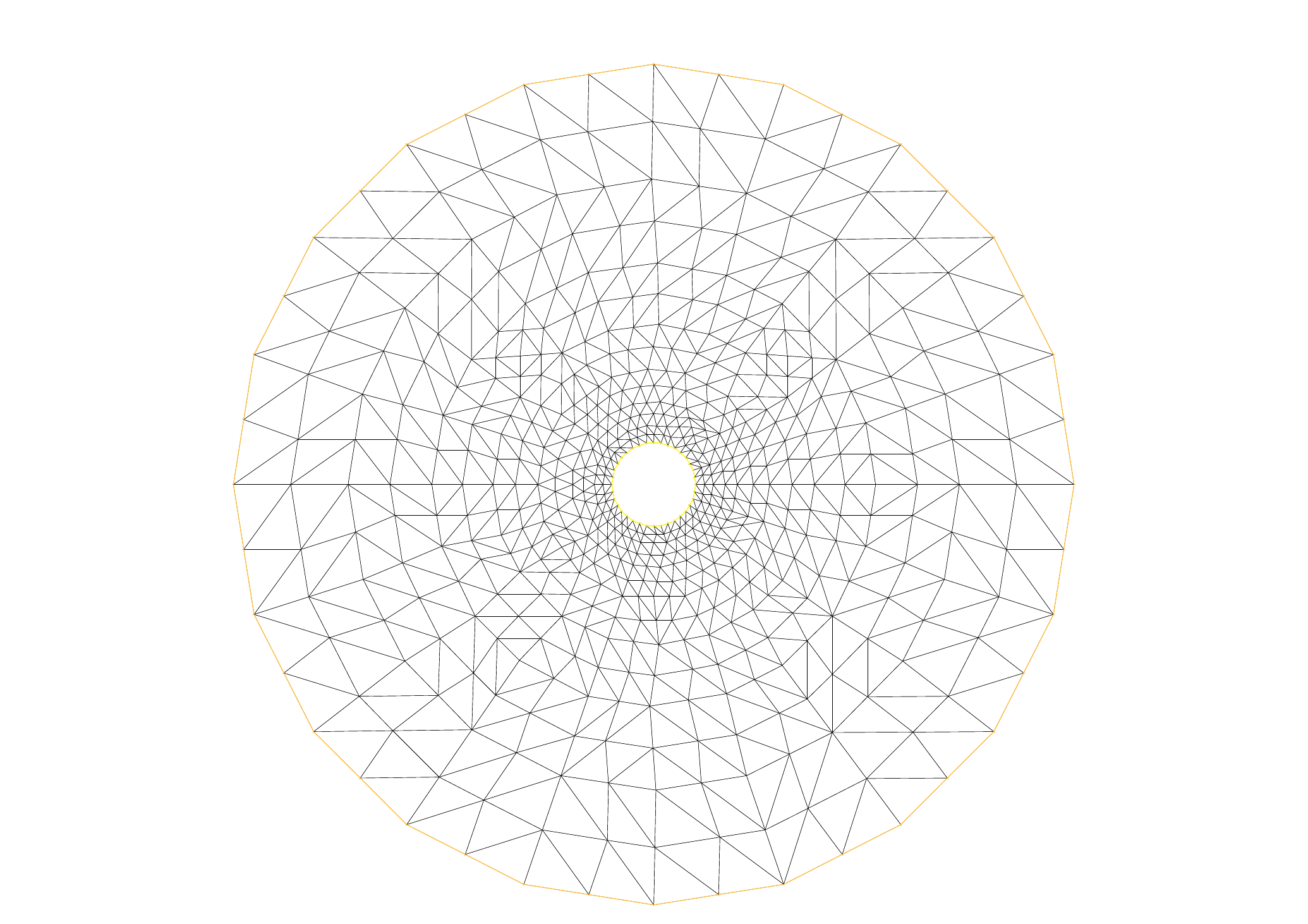}
    \caption{Mesh Level 2}
  \end{subfigure}
 
  \caption{Mesh for flow between an annulus}
  \label{fig:Grids}
\end{figure}

In this study, we have three relatively coarse meshes for the data, and one fine mesh for the DNS.     Our coarsest mesh, Mesh Level 1, is parameterized by 20  mesh points  around the outer circle and 18 mesh points  around the immersed circle, then extended to all of $\Omega$ as a Delaunay mesh \cite{DRSbook}. 
Mesh Level 2 is  generated  by splitting each triangle in Mesh Level 1 into $4$ sub-triangles.   Each grid in Mesh Level 2 is  then split into  $4$  sub-triangles to generate Mesh Level 4.   To create the finest  mesh, Mesh Level 8,  all triangles in Mesh Level 4 are split into $4$.  Figure \ref{fig:Grids} shows the first two levels.   More details about the grids can be found in Table \ref{table:1}. 
In what follows, $h$ stands for the finest mesh (Mesh Level 8) and $H$ represents a coarse spatial resolution where the observational data are collected.

\begin{table}[h!]
\centering
\begin{tabular}{||c c c c||} 
 \hline
 Mesh's Type & $\#$ Vertices &  $\#$ Triangles & Max size of mesh \\ [0.5ex] 
 \hline\hline
 Mesh Level 1 & 164 & 290 & 0.389  \\ 
Mesh Level 2 & 618 & 1160  & 0.194 \\
Mesh Level 4 & 2396 & 4640 &  0.097 \\
Mesh Level 8 & 9432 & 18560  &  0.048 \\ [1ex] 
 \hline
\end{tabular}
\caption{Details on the spatial meshes}
\label{table:1}
\end{table}

\subsection{Reference Solution}
 Since we do not have access to a true solution for this problem, we use the  DNS  solution as our reference solution.  We first resolve all scales  down to Kolmogorov dissipation scales  $h \leq \Rey^{- \frac{1}{2}}$, Mesh Level 8,  and then compute the DNS solution $\bfu^h_{n+1}$ using  the first order scheme Algorithm  \ref{Algorithm} but without nudging  ($\mu = 0$).  
 The Taylor-Hood mixed finite elements are utilized for discretization in space on a Delaunay-Vornoi generated triangular mesh.
 The DNS is run from $t_0=-5$ to a final time $T= 100 $ with time step $\Delta t = 0.01$ starting from rest, i.e., $\bfu^h(t_0) =\bfzero$.
 
\begin{figure}[h!]
  \centering
  \begin{subfigure}[b]{0.4\linewidth}
    \includegraphics[width=\linewidth]{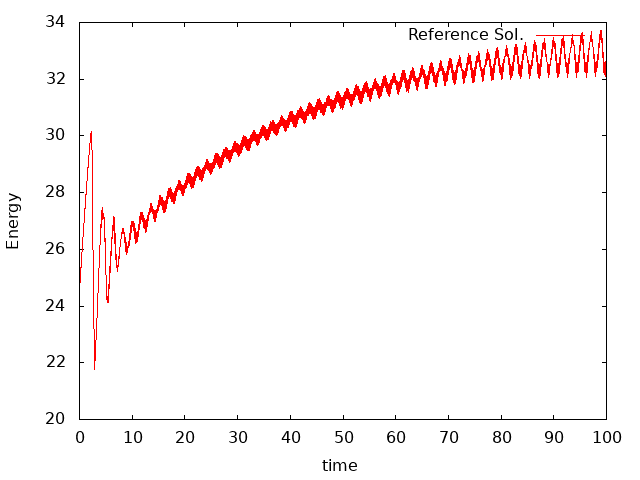}
  \end{subfigure}
 \begin{subfigure}[b]{0.4\linewidth}
    \includegraphics[width=\linewidth]{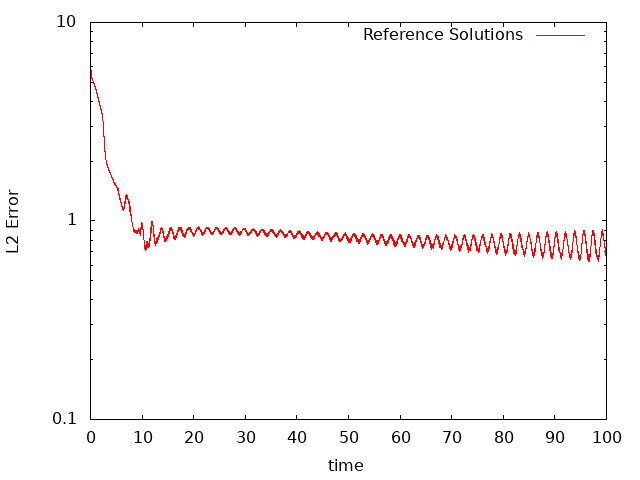}
 \end{subfigure}
 \caption{Left:  kinetic energy of the reference solution $\bfu^h$. Right: $\|\bfu_1^h - \bfu_2^h\|_{L^2}$;  $\bfu^h_1$  and $\bfu^h_2$ two reference solutions with different initial conditions}
   \label{fig:Shear-True}
   \end{figure}

The kinetic energy time series in Figure \ref{fig:Shear-True} shows the solution settling into regular oscillations.   Without the true initial condition, then, one can expect at least a lag in another solution compared to the reference solution as shown on the right.  More over,  the $L^2$ norm of difference of two DNS solutions  with two different initial conditions in Figure  \ref{fig:Shear-True} indicates the sensitivity of the solution to the initial conditions as shown on the right.

\subsection{Nudged Solution}
To compute the solution to \eqref{IMEX-FEM}, \eqref{incomp}, we start from zero initial conditions, i.e., $\bfv^h_0 = \bfzero$,  with  $\mu = 100$  use the same spatial and temporal discretization parameters as for the DNS, and start nudging with the DNS solution.  Interpolation $I_H$  is carried out on the  different refinement levels of spatial grids, while the equations \eqref{IMEX-FEM}, \eqref{incomp} are solved on the finest mesh. 

To compare the effect of the higher order interpolation on the approximate solution, we consider linear and quadratic Lagrange interpolation.    For simplicity,  we first describe the idea locally on triangle with  observational data  available  at six  nodes; three at the vertices and  three at midpoint of the edges of the triangle.   For the interpolant, there are two options as shown in figure \ref{Fig;Interpolation}
\begin{enumerate}
\item Quadratic interpolation using the six nodes, 
\item Refining the  triangle to four sub-triangles. Linear interpolation on each of the four sub-triangles. 
\begin{figure}[h!]
  \centering
  \begin{subfigure}[b]{0.25\linewidth}
    \includegraphics[width=\linewidth]{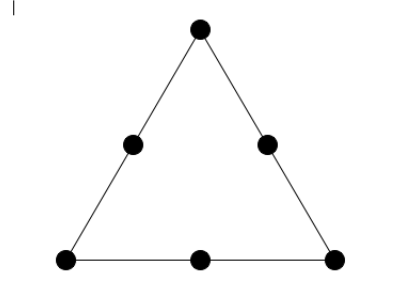}
  \end{subfigure}
  \begin{subfigure}[b]{0.25\linewidth}
    \includegraphics[width=\linewidth]{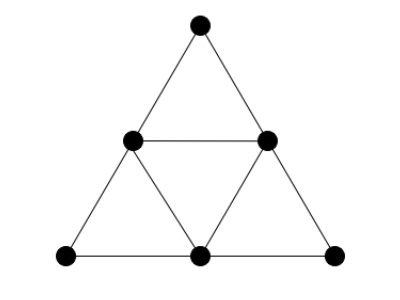}
   \end{subfigure}
   \caption{Quadratic  interpolation using six nodes (left) versus linear interpolation on finer grid (right)}
    \label{Fig;Interpolation}
   \end{figure}
\end{enumerate}

\begin{figure}[h!]
  \centering
  \begin{subfigure}[b]{0.4\linewidth}
    \includegraphics[width=\linewidth]{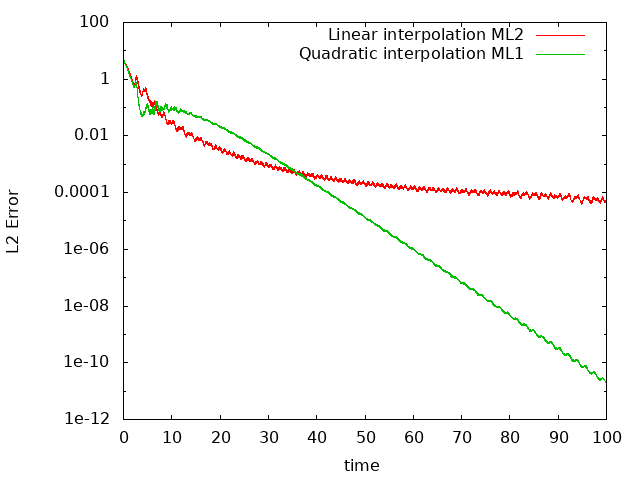}
  \end{subfigure}
  \begin{subfigure}[b]{0.4\linewidth}
    \includegraphics[width=\linewidth]{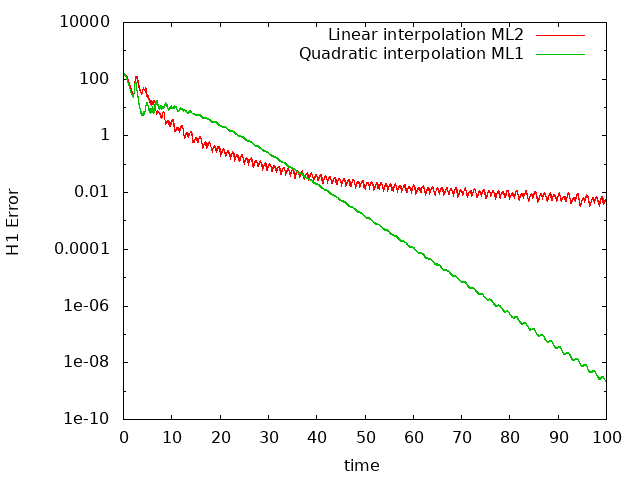}
   \end{subfigure}
   \caption{Data in Mesh Level 2. $L^2$ (left) $H^1$ (right) norm Error; Linear interpolation on Mesh Level 2 Vs. Quadratic interpolation on Mesh Level 1}
   \label{fig:ShearML12}
   \end{figure}

\begin{figure}
   \centerline{\includegraphics[scale=.45]{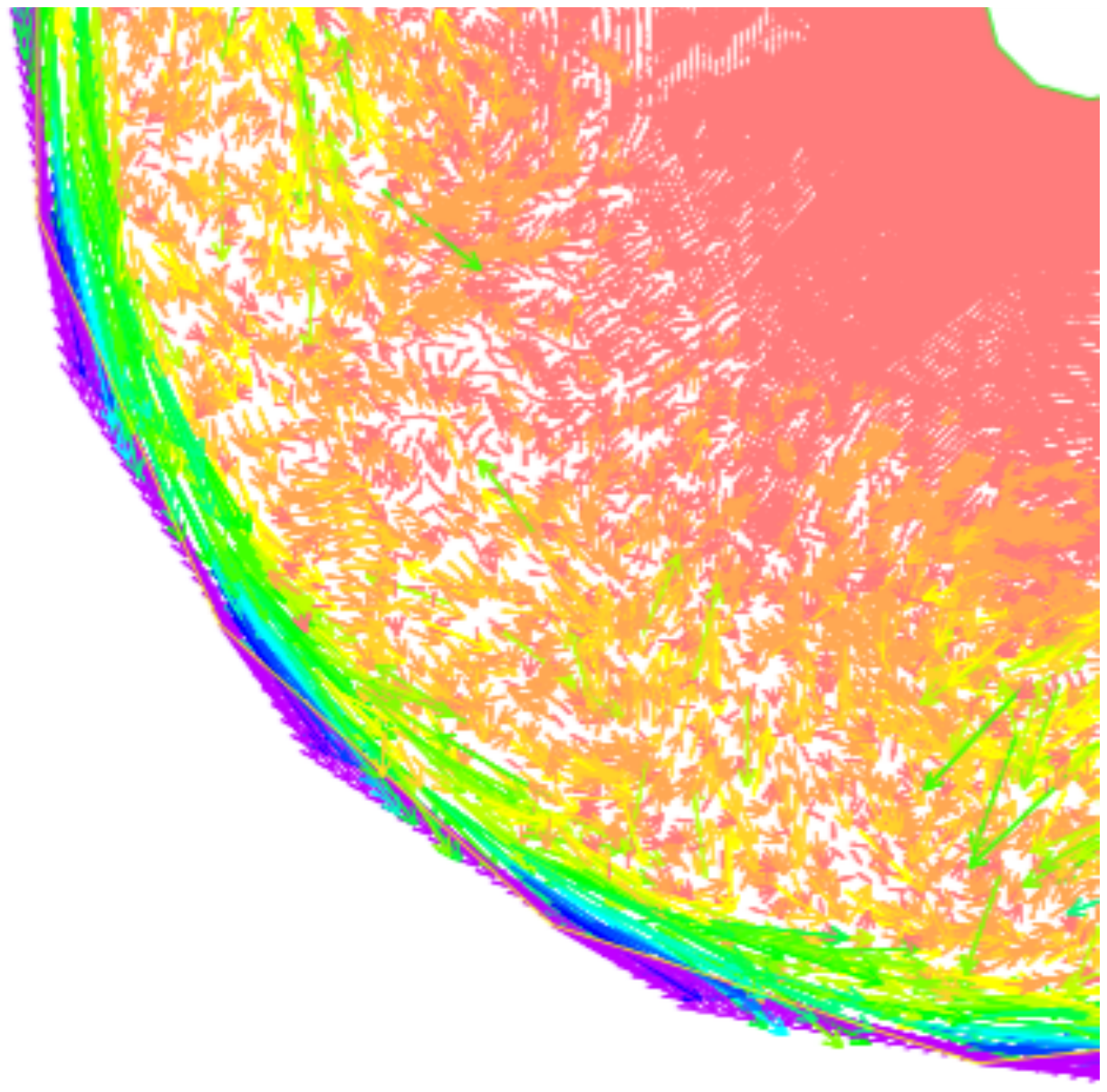} \hskip -1.35truein \includegraphics[scale=.45]{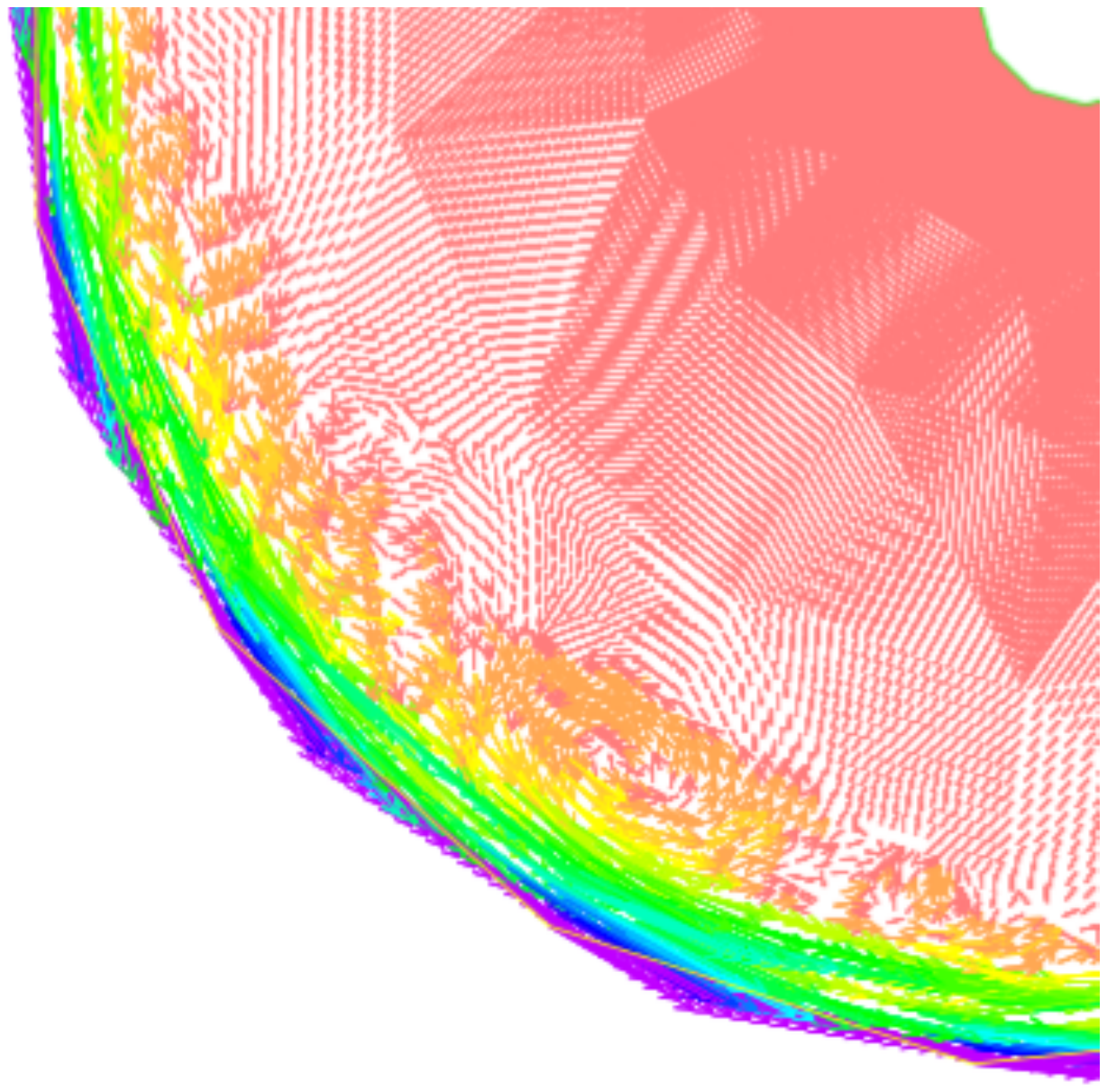} \hskip -1.35truein \includegraphics[scale=.45]{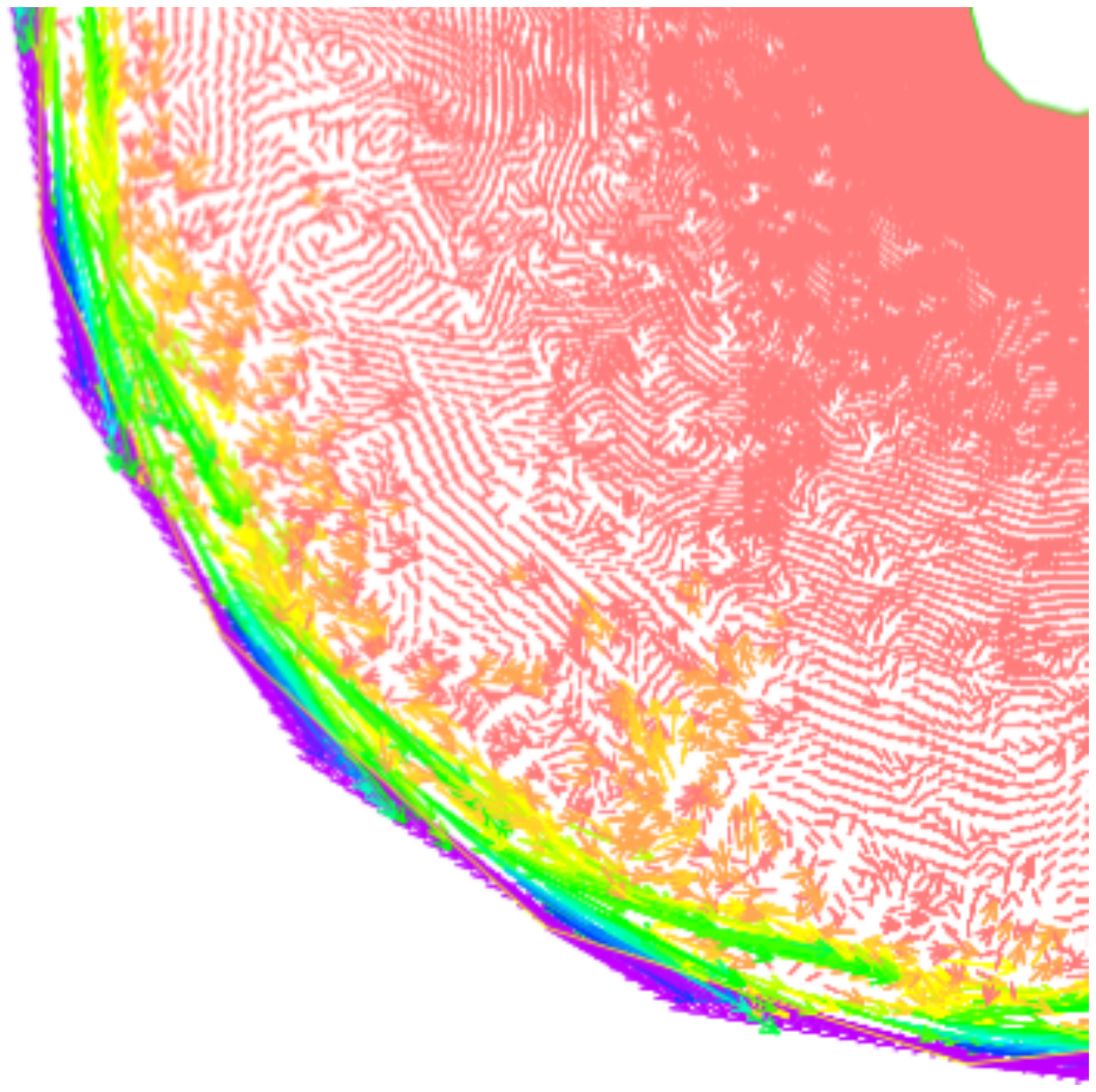}}
\vskip -1truein
\caption{Velocity vector, shear flow at $t=5$. Left: P1 interp. Mesh Level 2. Center: DNS (reference solution).  Right: P2 interp. Mesh Level 1. Color indicates vector length: watermelon=0, purple=12}
\label{shearfig}
\end{figure}

\begin{figure}[h!]
  \centering
  \begin{subfigure}[b]{0.4\linewidth}
    \includegraphics[width=\linewidth]{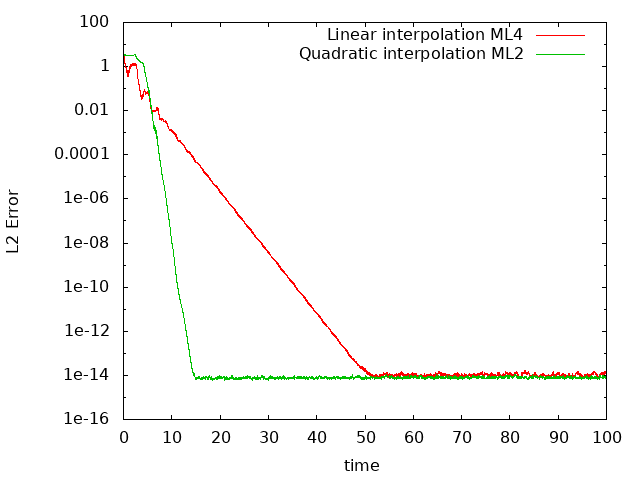}
  \end{subfigure}
  \begin{subfigure}[b]{0.4\linewidth}
    \includegraphics[width=\linewidth]{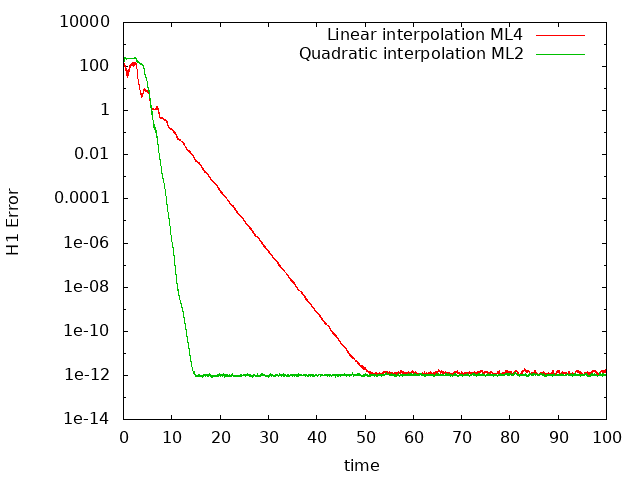}
    \end{subfigure}
   \caption{Data in Mesh level 4. $L^2$ (left) $H^1$ (right) norm Error; Linear interpolation on Mesh Level 4 Vs. Quadratic interpolation on Mesh Level 2}
   \label{fig:ShearML24}
   \end{figure}

The time evolution of the $L^2$ and $H^1$- norms for the difference in velocities of the reference and nudged solutions are shown in Figures \ref{fig:ShearML12}, and \ref{fig:ShearML24}.  The simulations are made for linear and quadratic Lagrange interpolate with two data resolutions to determine the effect of the higher order interpolation of the error.  
In each case the interpolants use the same set of observed data. In short, in all  mesh refinements,  the synchronization is made at a better rate using quadratic interpolation in compare with the linear one. For data in Mesh Level 2, nudging with the quadratic interpolant achieves synchronization at an exponential rate, while the error using linear interpolation decays much slower, Figure \ref{fig:ShearML12}.  Figure \ref{shearfig} shows a zoom of the velocity vectors in the lower left quadrant of the domain, comparing the nudging approaches on Mesh level 2 with the result of DNS.  The other quadrants are similar.  Both methods  synchronize to machine precision using data on Mesh Level 4, but the quadratic interpolant does so in roughly one-third the time, Figure \ref{fig:ShearML24}.   For data on Mesh Level 8, i.e., the irrelevant case of full knowledge of the reference solution, the errors are nearly the same, hence not shown here.



\newpage
\section{Computational Study \RNum{2} ; Body-forced case}

In this study, we consider the two-dimensional flow between two offset circles. The domain is a disk with  the omission of a smaller off-center disc inside given by
$$
\Omega = \{ (x,y) \in \mathbb{R}^2: x^2 + y^2 \leq 1^2 \hspace{0.3cm} \text{and}  \hspace{0.3cm} (x- 0.5)^2 + y^2 \geq (0.1)^2 \}.
$$
No-slip, no-penetration boundary conditions are imposed on both circles, and  the  flow is driven by the counterclockwise rotational  time-independent body force
$$\bff (x,y) = (- 4y (1- x^2 - y^2), 4 x (1-x^2-y^2))^{\tran}.$$

As the
flow rotates about  the origin it interacts with the inner boundary generating complex flow structures including the
formation of a Von Karman vortex street. This vortex street rotates and itself re-interacts with the immersed circle, creating more complex structures.  All the simulations are run over the time interval $[0, 40]$ with a time-step size $\Delta t = 0.01$ and Reynolds number $ \Rey = 600$, the same as done  in \cite{JL14}.  We use the same mesh structures generated in the last section adapted to the new domain; one fine mesh for DNS, and three relatively coarse meshes levels for data.

\subsection{Reference Solution} Again, since we do not have access to a true solution for this problem, we instead run a DNS and use the computed solution $\bfu^h$ as the reference solution.   
The  initial condition $\bfu_0$, 
is generated by solving the steady Stokes problem with  the same body forces $\bff (x,y)$.   

\begin{figure}[h!]
  \centering
  \begin{subfigure}[b]{0.4\linewidth}
    \includegraphics[width=\linewidth]{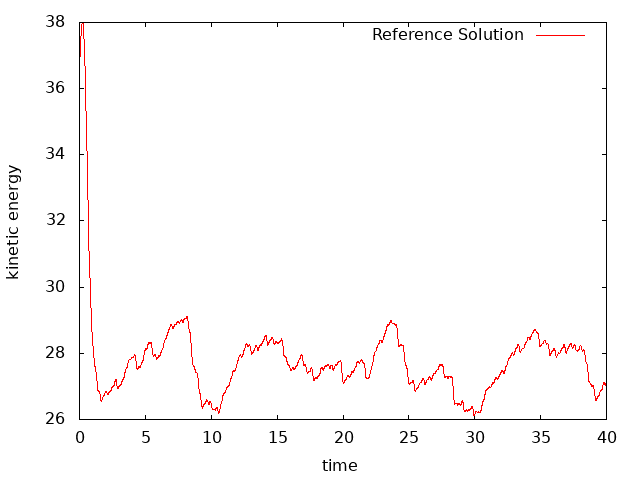}
  \end{subfigure}
  \begin{subfigure}[b]{0.4\linewidth}
    \includegraphics[width=\linewidth]{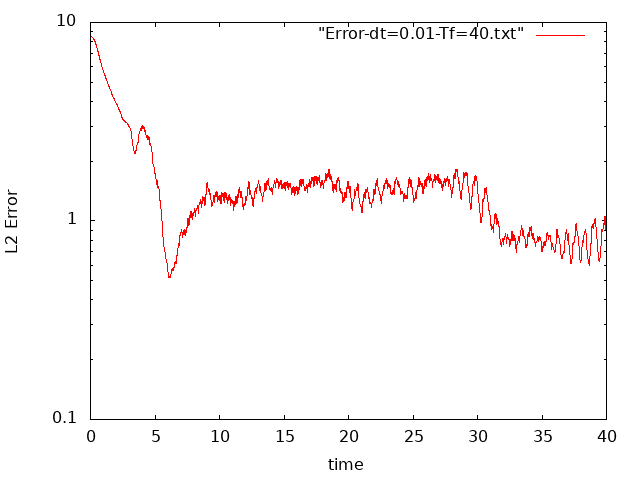}
   \end{subfigure}
   \caption{Left: Kinetic Energy of the reference  solution, indicating chaos. Right: $\|\bfu_1^h - {\bfu}_2^h\|_{L^2}$,  $\bfu_1^h$ and  ${\bfu}_2^h$ two reference solutions with different initial conditions}
    \label{fig:True}
   \end{figure}

The fluctuations of the reference solution's  kinetic energy  in time, shown in Figure \ref{fig:True},  is  evidence of  chaotic behavior.   The sensitivity with respect to initial data is demonstrated by DNS with two different initial conditions, $\bfu_{0_{1}} =$ generated by solving the steady state Stokes and $\bfu_{0_{2}} = \textbf{0}$.  The difference of the two solutions is plotted in Figure \ref{fig:True} right. This underscores the significance of a nudged solution synchronizing with the reference solution.


\subsection{Nudged Solution}
\begin{figure}[h]
  \centerline{\includegraphics[scale=.45]{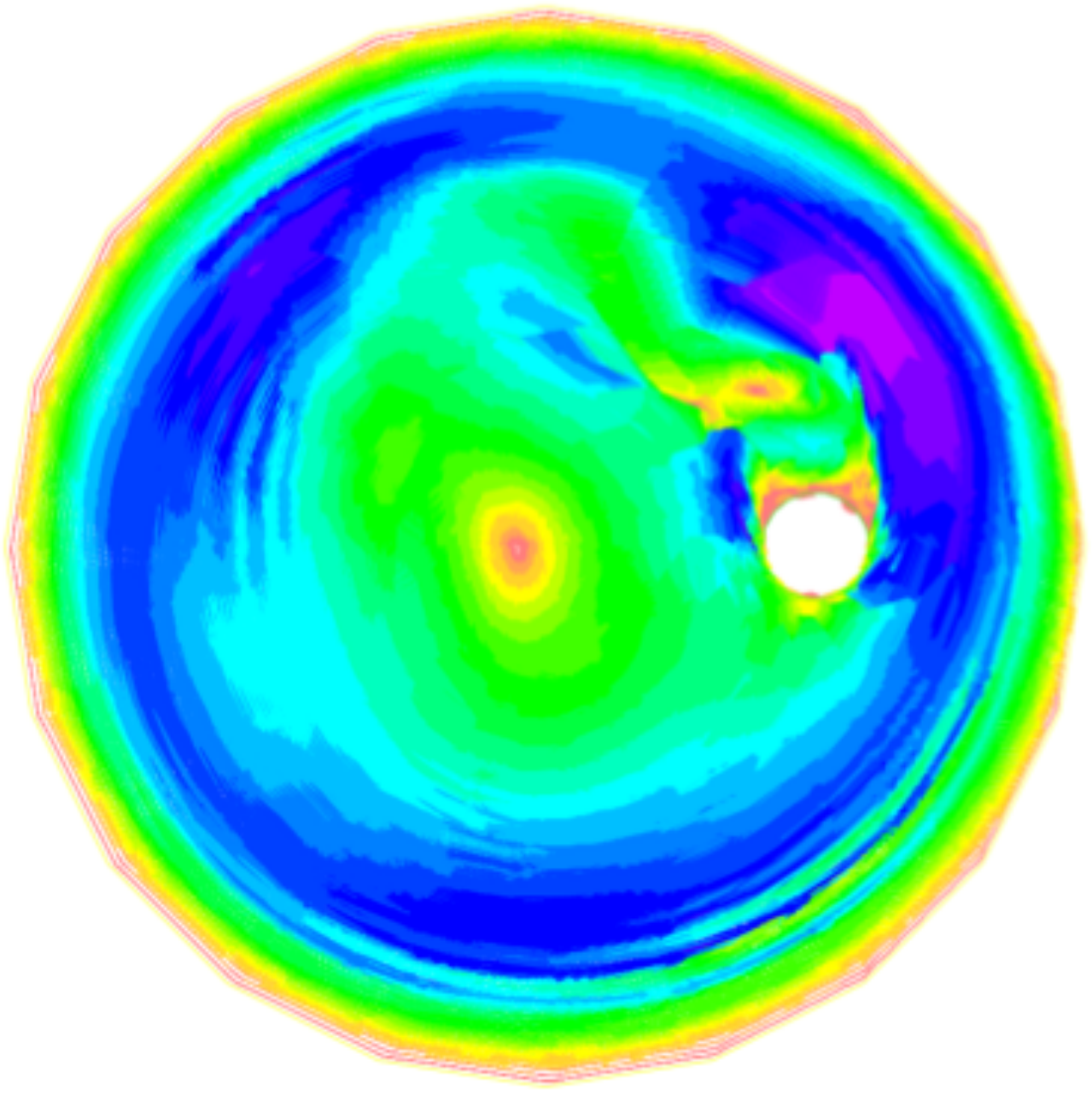} \hskip -1.5truein \includegraphics[scale=.45]{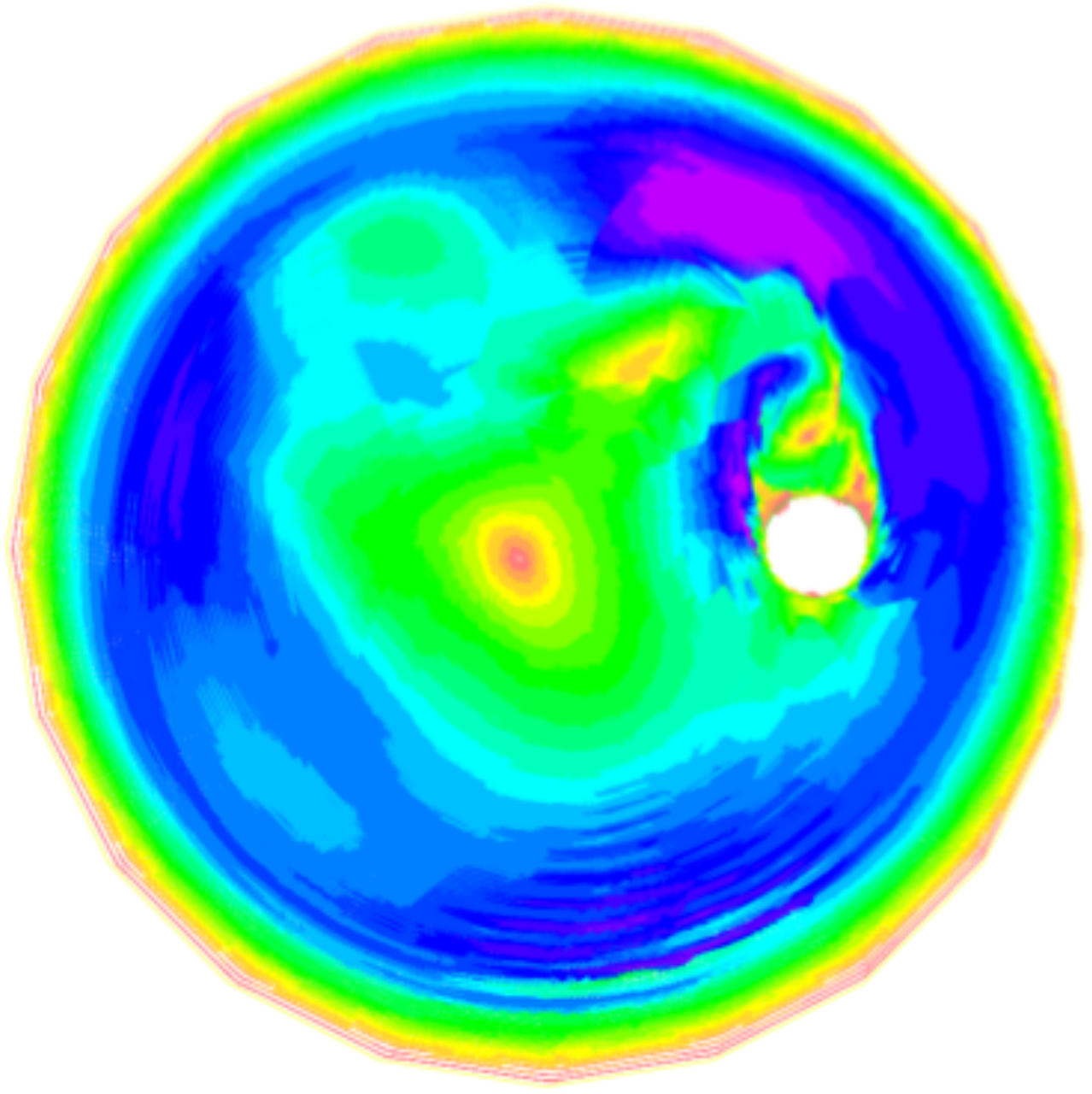} \hskip -1.5truein \includegraphics[scale=.45]{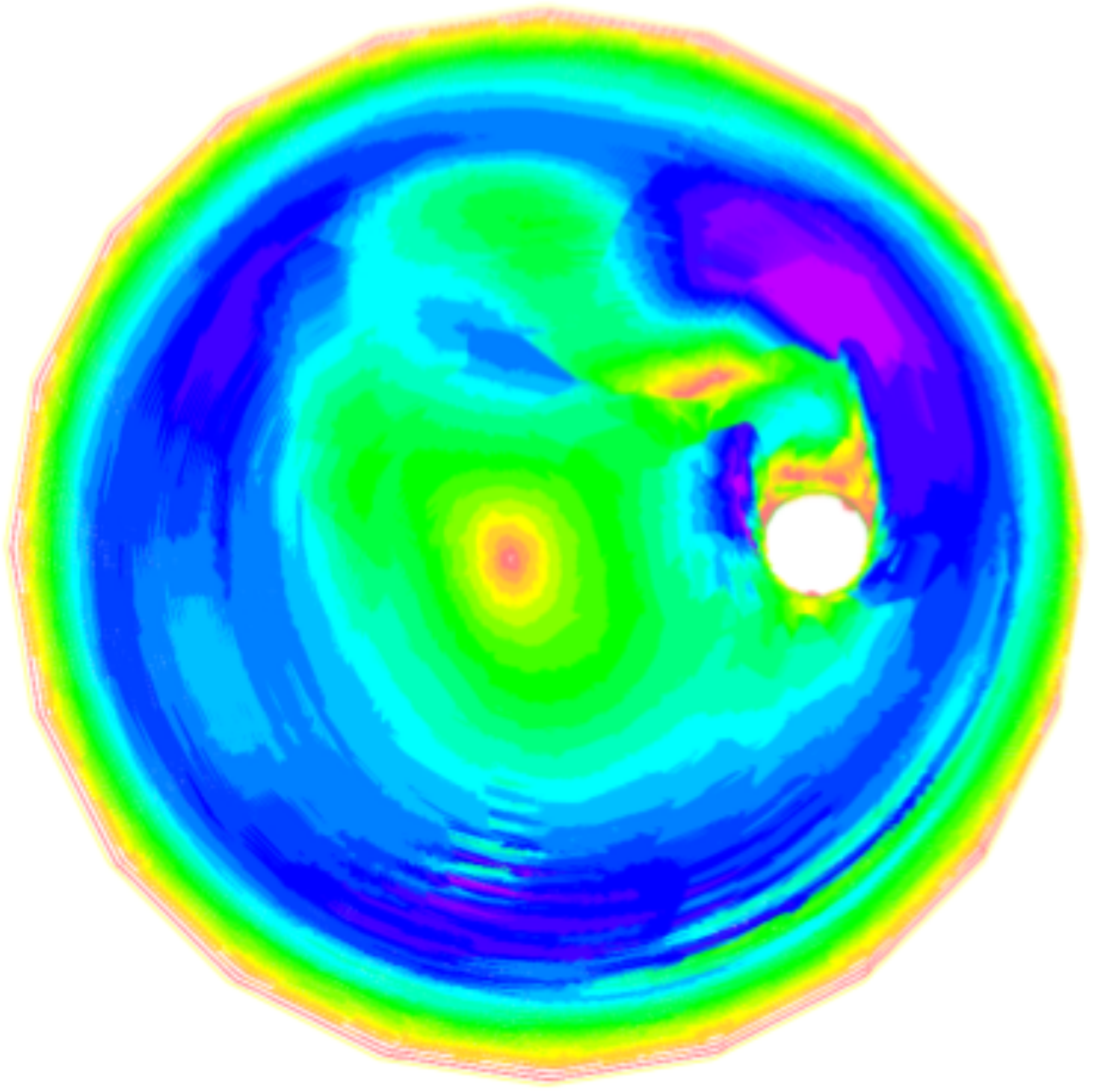}}
  \vskip -2.5truein\centerline{\includegraphics[scale=.45]{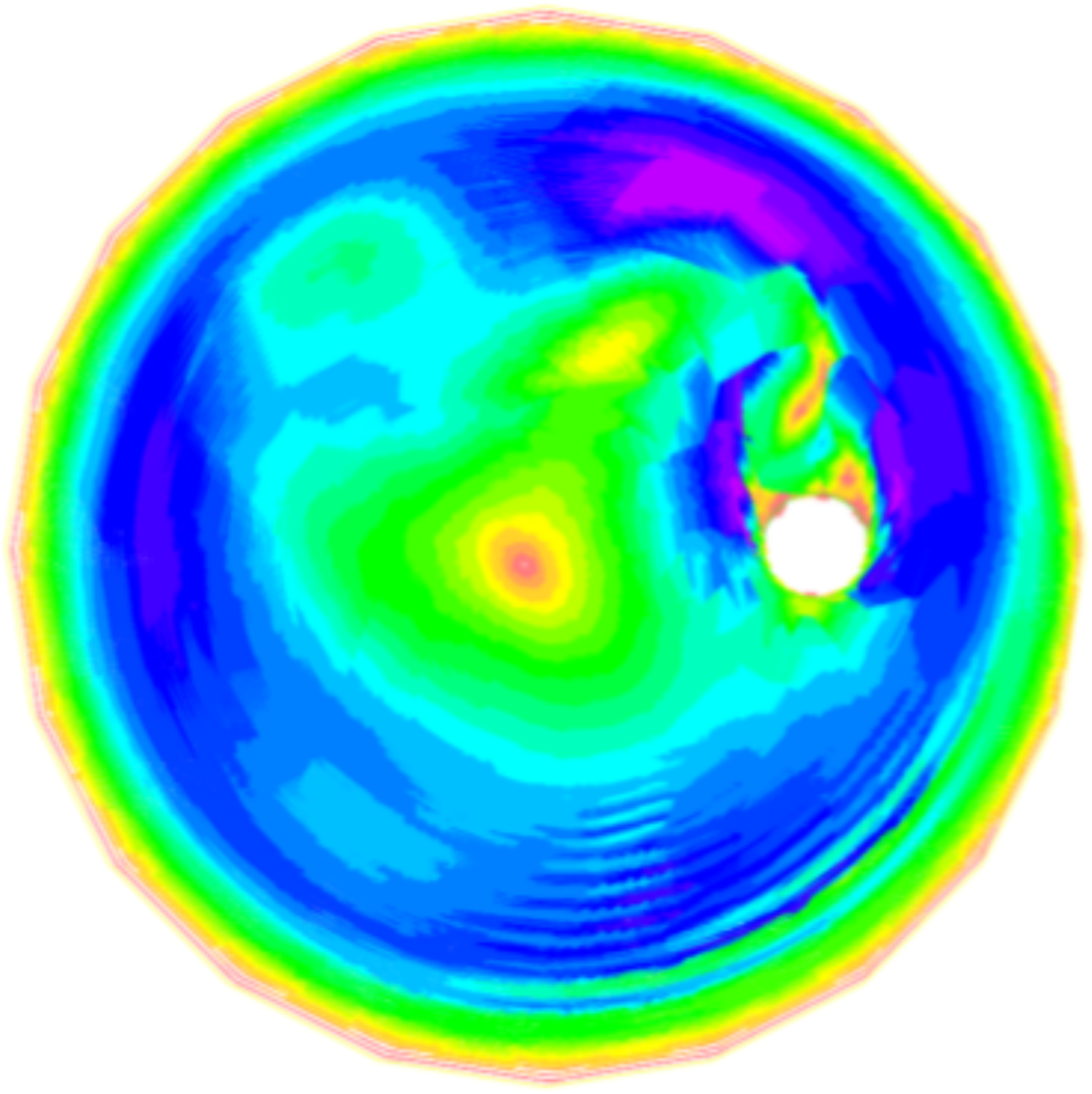} \hskip -1.5truein \includegraphics[scale=.45]{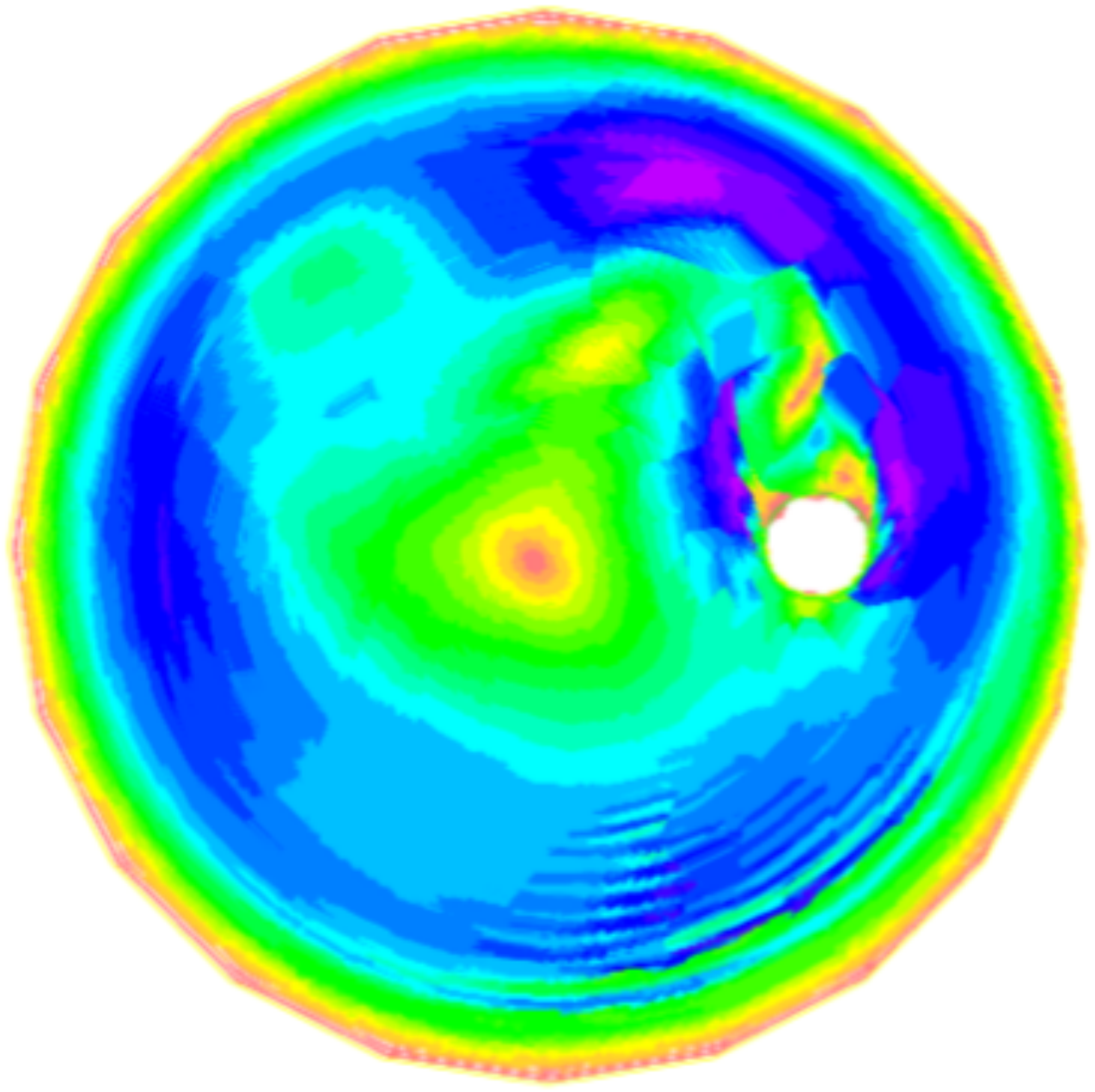} \hskip -1.5truein \includegraphics[scale=.45]{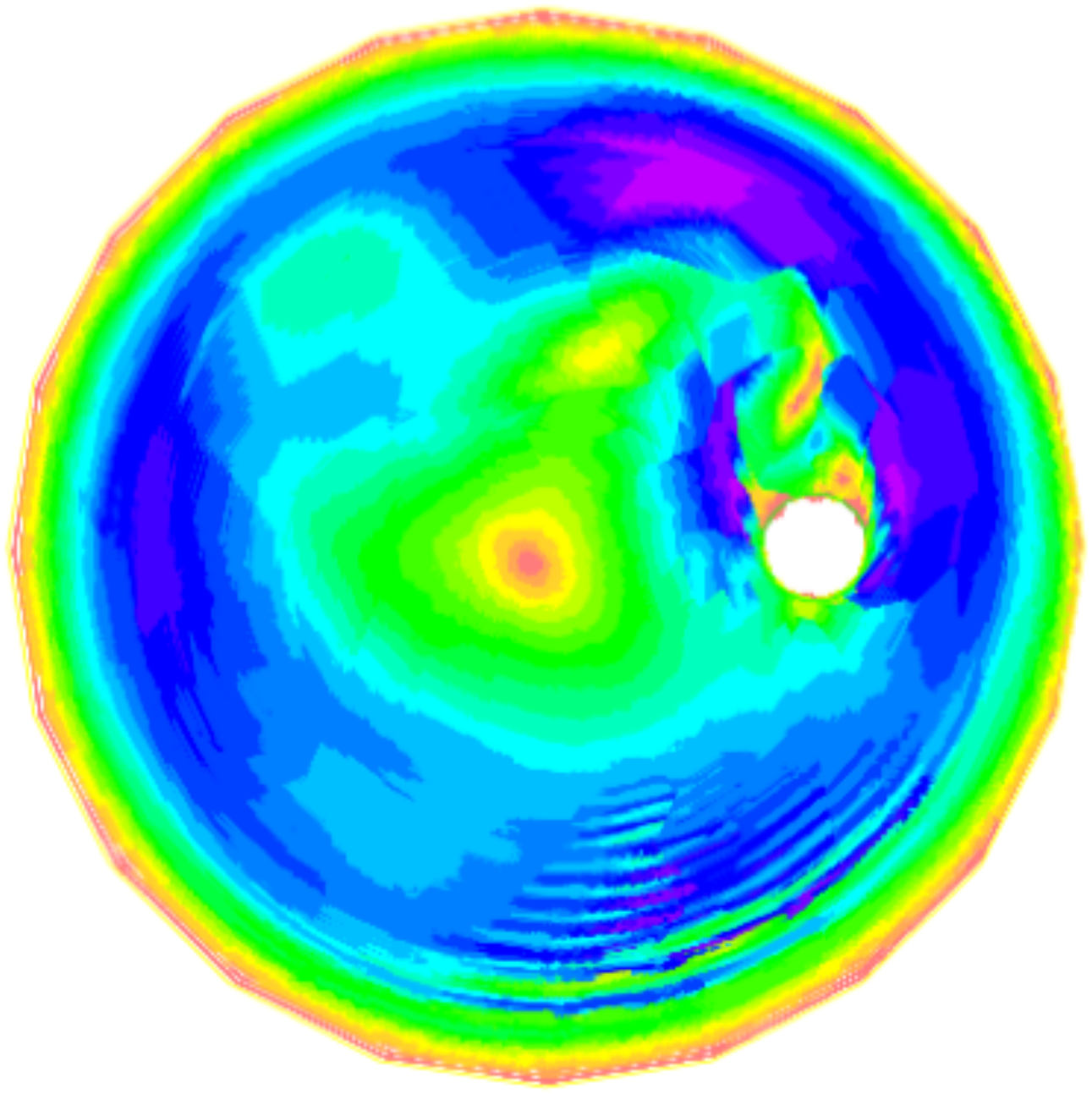}}
 \vskip -1.25truein 
\caption{Velocity vector, body force flow.  Top $t=5$, bottom $t=10$. Left: P1 interp. Mesh Level 2. Center: DNS (reference solution).  Right: P2 interp. Mesh Level1. Color indicates vector length: orange $=0$,  blue $\approx 5$, purple $\approx 7$ }
\label{forcefig}
\end{figure}

\begin{figure}[h!]
  \centering
  \begin{subfigure}[b]{0.4\linewidth}
    \includegraphics[width=\linewidth]{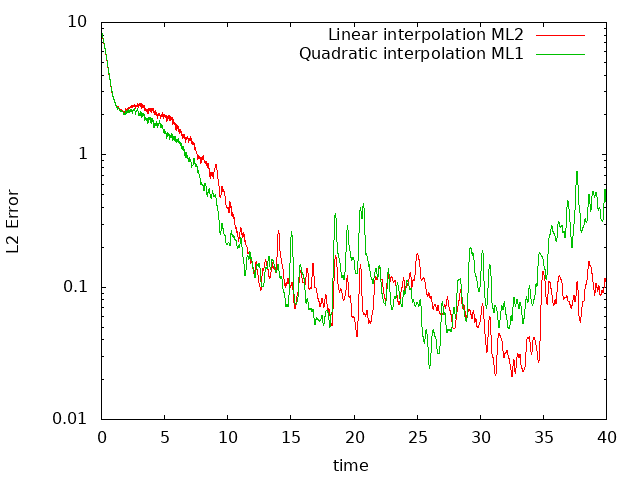}
  \end{subfigure}
  \begin{subfigure}[b]{0.4\linewidth}
    \includegraphics[width=\linewidth]{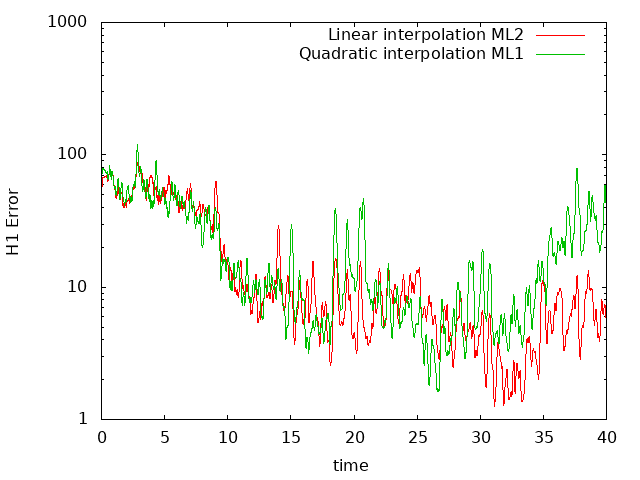}
   \end{subfigure}
   \caption{Data in Mesh Level 2.  $L^2$ (left) $H^1$ (right) norm Error; Linear interpolation on Mesh Level 2 Vs. Quadratic interpolation on Mesh Level 1}
   \label{fig:P1ML2-P2ML1}
   \end{figure}

\begin{figure}[h!]
  \centering
  \begin{subfigure}[b]{0.4\linewidth}
    \includegraphics[width=\linewidth]{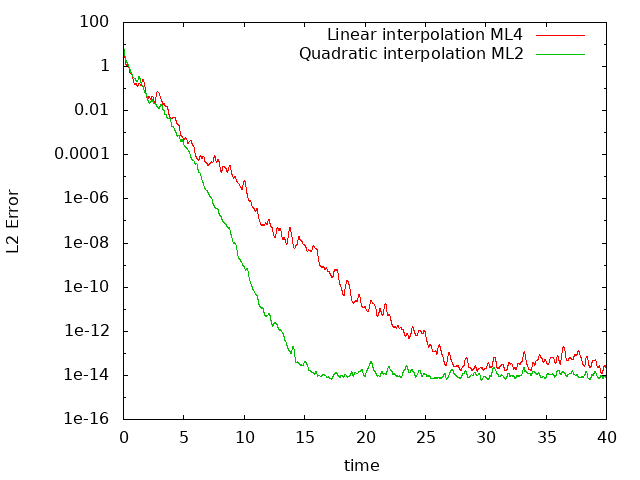}
  \end{subfigure}
  \begin{subfigure}[b]{0.4\linewidth}
    \includegraphics[width=\linewidth]{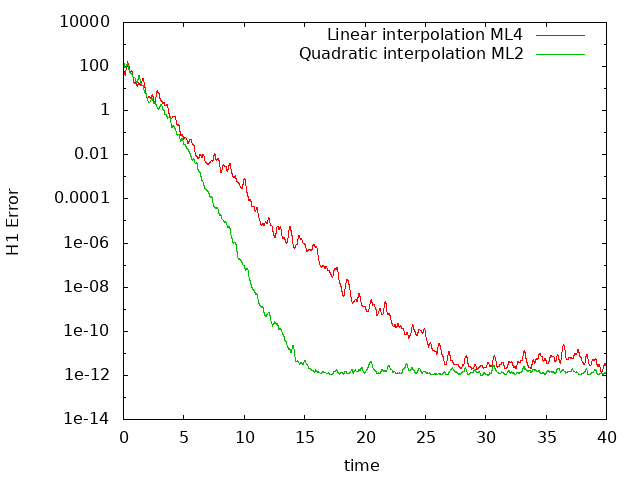}
    \end{subfigure}
   \caption{Data in Mesh Level 4.  $L^2$ (left) $H^1$ (right) norm Error; Linear interpolation on Mesh Level 4 Vs. Quadratic interpolation on Mesh Level 2}
   \label{fig:P1ML4-P2ML2}
   \end{figure}

   As in the shear flow case we start from zero initial conditions  $\bfv^h(0) = \textbf{0}$, set  $\mu =10$,  and use the same spatial and temporal discretization parameters as the DNS. While the nudging solution $\bfv^h$ is computed on the fine mesh, Mesh Level 8, the interpolation $I_H$  is either  linear on the Mesh Level $2k$  or  quadratic  on  the Mesh Level $k$, for  $k =1,2$. As before,  in both situations, we use the  same amount of observational data (i.e. locally six nodes on each triangle).  The $L^2$ and $H^1$ errors are plotted in Figures \ref{fig:P1ML2-P2ML1}, and \ref{fig:P1ML4-P2ML2}.  In this case the data is too sparse on Mesh Level 2, for interpolation of either degree to synchronize to near machine precision.  Figure \ref{forcefig} and \ref{fig:P1ML2-P2ML1}  shows that although the quadratic interpolation has slightly better performance at the beginning, nevertheless, both nudging captures the main features of the velocity field.   Both interpolation methods synchronize using data on Mesh Level 4, with quadratic interpolation doing so in about one-half the time, Figure \ref{fig:P1ML4-P2ML2} .  As expected and like the shear flow case,  errors are almost the same with complete knowledge of the flow,  not shown here.
   
   
\newpage

\section{Continuous time analysis}\label{SecAnalysis}

We have shown that quadratic interpolation can, in practice, outperform linear interpolation when used in nudging.   A complete analysis to support this would be one stated in terms of bounds on 
$\|\bfu-\bfv^h\|_0$, as is done for constant and linear interpolation in \cite{AOT14}.  To the best of our knowledge, these results have not yet been extended to higher order interpolation.  While doing so, i.e., establishing that $\|\bfu-\bfv^h\|_0 \to 0$ as $t\to \infty$, when nudging with higher order interpolation may be straightforward, deriving an estimate that indicates an advantage of higher order interpolation is another matter. In this section, we sketch a proof of synchronization using higher order interpolation in the simpler case of a semi-discrete, continuous in time finite element approximation, and at the same time, illustrate the difficulty in demonstrating such an advantage over lower order interpolation. 

To this end, we consider a reference solution satisfying the semi-discrete approximation
 \begin{equation}\label{NSE-FEM}
 \begin{split}
 (\bfu^h_t , \Theta^h) + \nu \, (\nabla  \bfu^h, \nabla  \Theta^h) + b(\bfu^h, \bfu^h, \Theta^h) = (f, \Theta^h) &,\\
 (\bfu^h (\cdot,0) - \bfu_0\,  , \, \Theta^h) = 0&,
 \end{split}
 \end{equation}
for all $\Theta^h \in \bf  V^h$, where, thanks to (\ref{LBBH}), the pressure has been eliminated by restricting the velocity to the space of discrete, divergence free functions $\bfV^h$.  The continuous time nudging scheme is given by
 \begin{equation}\label{DA-FEM}
 \begin{split}
 (\bfv^h_t , \Theta^h) + \nu \, (\nabla  \bfv^h, \nabla  \Theta^h) + b(\bfv^h, \bfv^h, \Theta^h) = (f, \Theta^h) - \mu \,  (I^k_H(\bfv^h - \bfu^h), \Theta^h)&,\\
 (\bfv^h (\cdot,0) - \bfv_0\,  , \, \Theta^h) = 0&,
 \end{split}
 \end{equation}
  for any $\Theta^h \in \bf  V^h$, where $I_H^k$ is $k^{\text{th}}$ order (degree)  interpolating polynomial satisfying (\ref{Ineq:Intrpolate}).    
  
With the following decomposition of the nonlinear term  
\begin{equation*}
\begin{split}
b(\bfu^h, \bfu^h, \Theta^h) - b(\bfv^h, \bfv^h, \Theta^h) & =b(\bfu^h, \bfu^h, \Theta^h) -  b(\bfv^h, \bfu^h, \Theta^h)\\
&+  b(\bfv^h, \bfu^h, \Theta^h) - b(\bfv^h, \bfv^h, \Theta^h)\\
&= b(\bfw^h, \bfu^h, \Theta^h) + b(\bfv^h, \bfw^h, \Theta^h),
\end{split}
\end{equation*}
the difference $\bfw^h = \bfu^h - \bfv^h \in \bfV^h$ satisfies
  $$(\bfw^h_t , \Theta^h) + \nu \, (\nabla  \bfw^h, \nabla  \Theta^h) + b(\bfw^h, \bfu^h, \Theta^h) + b(\bfv^h, \bfw^h, \Theta^h) =- \mu \,  (I_H \bfw^h , \Theta^h).$$
Setting $\Theta^h = \bfw^h$, the second nonlinear term vanishes because $b(\cdot , \cdot, \cdot)$ is explicitly skew-symmetrized. Thus 

\begin{equation} \label{Er-Eq1}
\frac{1}{2}\, \frac{d}{dt} \|\bfw^h\|_0^2 +  \nu\, \|\bfw^h\|^2_{1} =   - \mu (I_H^k \bfw^h , \bfw^h) -  b(\bfw^h, \bfu^h, \bfw^h).
\end{equation}
Using  Young's inequality,  along with \eqref{Ineq:Intrpolate} and  \eqref{Ineq:Inverse}, we have
\begin{equation} \label{Er-Eq2}
\begin{split}
- \mu (I_H^k \bfw^h , \bfw^h) & =  - \mu (I_H^k \bfw^h  - \bfw^h + \bfw^h , \bfw^h) =  \mu (  \bfw^h  - I_H^k \bfw^h , \bfw^h) - \mu \|\bfw^h\|_0^2\\
& \leq \frac{\mu}{2}  \|\bfw^h  - I_H^k \bfw^h \|_0^2 + \frac{\mu}{2}\|\bfw^h \|_0^2 - \mu  \|\bfw^h \|_0^2\\
& \leq \frac{\mu}{2}\,  C_{k+1,0}^2 H^{2k+2} \|\bfw^h \|^2_{{k+1}} -  \frac{\mu}{2}\|\bfw^h \|_0^2\\
& \leq \frac{\mu}{2}\, C_{k+1,0}^2\, H^{2k+2} \,  \Cinv_{k+1,1}^2\, h^{-2k}  \|\bfw^h \|^2_{{1}} -  \frac{\mu}{2}\|\bfw^h \|_0^2\\
& \leq   \frac{\nu}{2}    \|\bfw^h \|^2_{1} -  \frac{\mu}{2}\|\bfw^h \|_0^2\;,\\
\end{split}
\end{equation}
provided 
\begin{align}\label{kcond}
\nu \ge \mu \,  C_{k+1,0}^2\,  \,  \Cinv_{k+1,1}^2\, \left(\frac{H}{h}\right)^{2k}\, H^{2}.\, 
\end{align}
 The H\"older, Ladyzhenskaya \footnote{$\|\varphi\|_{L^4} \leq \|\varphi\|^{\frac{1}{2}} \, \|\nabla \varphi\|^{\frac{1}{2}}.$} 
 , Poincar\'e  \footnote{ $\|\varphi\|_0 \leq \frac{1}{\sqrt{\lambda_1}}\, \|\varphi\|_1,$ \hspace{0.3cm}  
 $\|\varphi\|_1 \leq \frac{1}{\sqrt{\lambda_1}} \, \|\varphi\|_2.$
}
 , and Agmon \footnote{$\|\varphi\|_{L^\infty} \leq \|\varphi\|_0^{\frac{1}{2}} \, \| \varphi\|_2^{\frac{1}{2}}.$}
 inequalities give us the following estimate on the nonlinear term in (\ref{Er-Eq1})
\begin{equation} \label{Er-Eq4}
\begin{split}
& b(\bfw^h, \bfu^h, \bfw^h) = \frac{1}{2} (\bfw^h \cdot \nabla \bfu^h , \bfw^h) - \frac{1}{2} (\bfw^h \cdot \nabla \bfw^h , \bfu^h) \\
& \leq \frac{C}{2} \| \bfu^h\|_1 \|\bfw^h\|^2_{L^4} +  \frac{C}{2} \|\bfu^h\|_{L^{\infty}} \|\bfw^h\|_0 \|\bfw^h\|_1\\
&\leq \frac{C}{2} \| \bfu^h\|_1  \|\bfw^h\|_1  \| \bfw^h\|_0 + \frac{C}{2} \| \bfu^h\|_0^{\frac{1}{2}} \| \bfu^h\|^{\frac{1}{2}}_{2} \|\bfw^h\|_0 \|\bfw^h\|_1\\
& \leq \frac{\nu}{4}  \|\bfw^h\|^2_1  + \frac{C}{4 \nu} \| \bfu^h\|^2_1  \|\bfw^h\|_0^2  \, + \,  \frac{\nu}{4}  \|\bfw^h\|^2_1  + \frac{C}{4\nu} \| \bfu^h\|_0 \| \bfu^h\|_{2}\,  \|\bfw^h\|_0^2\\
& \leq \frac{\nu}{2}  \|\bfw^h\|^2_1  + \frac{C}{4 \nu} \left(\| \bfu^h\|^2_1  +  \| \bfu^h\|_0 \| \bfu^h\|_{2}\right) \|\bfw^h\|_0^2\\
& \leq \frac{\nu}{2}  \|\bfw^h\|^2_1  +   \frac{C}{2\nu \, \lambda_1}  \|\bfu^h\|_{2}^2\,  \|\bfw^h\|_0^2.
\end{split}
\end{equation}
Substituting \eqref{Er-Eq2}, \eqref{Er-Eq4} in 
\eqref{Er-Eq1}, we have
\begin{equation} \label{Er-Eq5}
 \frac{d}{dt} \|\bfw^h\|_0^2  + \left(\mu -  \frac{C}{\nu \lambda_1} \, \|\bfu^h\|_{2}^2 \right)  \|\bfw^h\|_0^2 \leq 0. 
\end{equation}

We next use the following uniform Gr\"onwall inequality proved in \cite{JonesTiti}.
 \begin{lem}\label{Gronwall} Let $\tau >0$ be arbitrary  and fixed. Suppose that $Y(t)$ is an
absolutely continuous function which is locally integrable such that 
$$\frac{d Y}{dt} + \alpha(t) Y \leq 0,$$
where 
$$\limsup_{t \rightarrow \infty}\, \int_t^{t+\tau} \alpha(s) ds \geq \gamma >0.$$
Then $Y(t) \rightarrow 0$ exponentially fast, as $t \rightarrow \infty.$
 \end{lem}
 
One can adapt an argument in \cite{FLT} for the solution to the NSE, to show that there exists a time $t_0 > 0$ such that for all $t \geq t_0$  and $\tau >0$
 \begin{align}\label{intH2bound}
 \int_t^{t+\tau} \|\bfu^h\|_2^2\, ds \leq  \big(c_0 e^{G^4}  + \tau \nu \lambda_1\big) \, \nu \lambda_1 G^2,
\end{align}
where $c_0$ is a positive non-dimensional constant. Now with $\alpha(s) = \mu -  \frac{C}{\nu \lambda_1} \, \|\bfu^h(s)\|_{2}^2$ in \eqref{Er-Eq5}, take  $\tau=\frac{1}{\nu\lambda}$, and assume

\begin{align}\label{mucond}\mu \ge 2\,  
 \big(c_0 e^{G^4}  + 1\big) \, \nu \lambda_1\, G^2,
\end{align} so that 
$$\int_t^{t+(\nu\lambda)^{-1}}\left( \mu -  \frac{C}{\nu \lambda_1} \, \|\bfu^h(s)\|_{2}^2\right) \, ds  \geq \big(c_0 e^{G^4}  + 1\big) \,  G^2 = \gamma >0,$$
and by Lemma \ref{Gronwall} it follows that $\|\bfw^h\|_0 \rightarrow 0$ at an exponential rate.

\begin{re}

Using an alternative approach, we can avoid the exponential factor in $G$, at the cost of cancelling $H^{-2}$ factor in \eqref{bothcond}.  Taking $\Theta=\bfu^h$ in \eqref{NSE-FEM}, one finds there exists $t_0>0$, such that for all $t>t_0$, the following well-known bound (see, e.g., (9.9) in \cite{Constantin-Navier1988}) holds \begin{equation}
\int_t^{t+ \tau} \|\bfu^h(s)\|_1 \ ds \leq \, \left( 2 + \tau \nu \lambda_1\right)\, \nu\,  G^2 \;.
\end{equation}
We can then apply the inverse inequality  \eqref{Ineq:Inverse}  to find
\begin{align}\label{intH2/1bound}
\int_t^{t+(\nu\lambda_1)^{-1}} \|\bfu^h\|_2^2\, ds \leq  \Cinv_{2,1}^2 h^{-2}  \int_t^{t+(\nu\lambda_1)^{-1}} \|\bfu^h\|_1^2\,\ ds \leq   3\Cinv_{2,1}^2 h^{-2}  \nu\, G^2\;.
\end{align}
Thus, if
\begin{align}\label{mucond2}\mu \ge 4 \, C \,  \Cinv_{2,1}^2\,  h^{-2}\, \nu \, G^2\;,
\end{align}
we have
$$\int_t^{t+(\nu\lambda)^{-1}} \left(\mu -  \frac{C}{\nu \lambda_1} \, \|\bfu^h(s)\|_{2}^2\right)\, ds  \geq  \frac{ 1}{ h^2\, \lambda_1} \,  C\, \Cinv^2_{2,1} \, G^2 = \gamma >0,$$
resulting in a range
\begin{align}\label{bothcond}
 4 \, C \,  \Cinv_{2,1}^2\,  h^{-2}\, \nu \, G^2 \le \mu \le  \nu\,  C_{k+1,0}^{-2} \, \Cinv_{k+1,1}^{-2}\left(\frac{h}{H}\right)^{2k} H^{-2}.
\end{align}
which, due to the factor of $h^{-2}$ in the lower bound, is clearly not achievable.  
\end{re}

\begin{re}
Similar analysis can be done to show synchronization in the $\bfH^1$ norm, by taking $\Theta^h=\Delta \bfw^h$.  In particular for $k=1$, one can find
$$
-\mu(I_H^1\bfw^h,\Delta \bfw^h) \le \frac{C_{1,0}^2\Cinv_{3,2}^2}{2 \nu} \mu \left(\frac{H}{h}\right)^2 \|\bfw^h\|_2^2 +
\frac{\nu}{2} \|\bfw^h\|_2^2 -\mu\|\bfw^h\|_1^2
$$
and exactly as in (3.22), (3.23) of \cite{FJT} 
\begin{align*}
|b(\bfu^h,\bfu^h,\Delta\bfw^h)-b(\bfv^h,\bfv^h,\Delta\bfw^h)| &\le |b(\bfw^h,\bfu^h,\Delta\bfw^h)|+|b(\bfv^h,\bfw^h,\Delta\bfw^h)| \\
& \le \frac{\nu}{4}\|\bfw^h\|_2^2 + \frac{C}{\nu}\left(1 + \frac{\|\bfu^h\|^2_2}{\lambda_1}\right) \|\bfw^h\|_1^2
\end{align*}
 One then can proceed as in the estimate of the $L^2$ error.
 \end{re}

 \section{Conclusion}

We have presented numerical simulations which demonstrate that nudging with higher order interpolation can synchronize when the data is too coarse for linear interpolation to do so, and does so faster when the data is fine enough for both to synchronize.  Two flows were tested: a shear flow in an annulus, and one with a body force in a disk with an off-center obstacle, both satisfying Dirichlet boundary conditions. We have shown rigorously that continuous data assimilation by nudging with higher order interpolation will synchronize with the solution of the spatially discretized NSE.  Even without estimating the error between the nudged solution and the actual solution of the NSE, the analysis is limited in gauging the benefit of using higher order interpolation.  The conditions \eqref{mucond} and \eqref{kcond} specify that a valid range for the relaxation parameter $\mu$ is
\begin{align}\label{bothcond2}
2\,  
 \big(c_0 e^{G^4}  + 1\big) \,  G^2 \le  \mu \le  \nu\,  C_{k+1,0}^{-2} \, \Cinv_{k+1,1}^{-2}\left(\frac{h}{H}\right)^{2k} H^{-2}\;.
\end{align}
Several points are in order.  
\begin{enumerate}
    \item  The condition \eqref{bothcond2} is achievable for order $k$ interpolation, provided $H$, the resolution of the data, is small enough.
    \item The higher the order $k$ is taken, the smaller $H$ is needed in the upper bound in  \eqref{bothcond}.
    This analysis is not sensitive enough to indicate the advantage of nudging with higher order interpolation demonstrated by our numerical experiments.
    \item The exponential factor in the Grashof number is due to assuming Dirichlet boundary conditions.  A much more reasonable, algebraic in $G$, lower bound would suffice for periodic boundary conditions. 
    Though in either case the restriction may seem impractical for turbulent flows that require large $G$, our simulations suggest it is far from sharp, as has been shown to be the case in comparing other numerical tests of nudging with corresponding analyses \cites{Cao-Algebraic2021,Farhat-Assimilation2018,Gesho-Acomputational2016,Hudson-Numerical2019}.  
\end{enumerate}

\end{document}